\documentclass{article}

\usepackage{amssymb}
\usepackage{amsmath}
\usepackage{color}
\usepackage{graphicx}
\usepackage{graphics}
\newtheorem{thm}{Theorem}
\newtheorem{prop}{Proposition}

\newtheorem{lem}{Lemma}
\newtheorem{rem}{Remark}

\newtheorem{defn}{Definition}

\newcommand{\RR}{{\mathbb R}}

\newcommand{\norm}[1]{\lVert#1\rVert}

\newcommand{\tr}[1]{{\text{Tr}(#1)}}

\newcommand{\dotex}{{\frac{d}{dt}}}

\begin{document}

\title{The geometry of low-rank Kalman filters}
\author{Silvere Bonnabel\thanks{S. Bonnabel is with Centre de Robotique, Math\'{e}matiques et Syst\`{e}mes,
      Mines ParisTech, Boulevard Saint-Michel 60, 75272 Paris, France.
      {\tt\small silvere.bonnabel@mines-paristech.fr}}
     \and Rodolphe Sepulchre\thanks{ Rodolphe Sepulchre is with Departement of Electrical Engineering and Computer Science,
      University of Li\`{e}ge, 4000 Li\`{e}ge, Belgium.
      {\tt\small \{r.sepulchre\}@ulg.ac.be}}}

\maketitle

\begin{abstract}
An important property of the Kalman filter is that the underlying Riccati flow is a contraction for the natural metric of the cone of symmetric positive definite matrices. The present paper studies the geometry of a low-rank version of the Kalman filter. The underlying Riccati flow evolves on the manifold of fixed rank symmetric positive semidefinite matrices. Contraction properties of the low-rank flow are studied by means of  a suitable metric recently introduced by the authors.
\end{abstract}

\section{Introduction}

The Kalman-Bucy filter \cite{Kalman-1961} (KF)  is a very popular method  in engineering, that allows  to
 compute an estimate of the state of a dynamical system from several
sensors measurements, possibly corrupted by measurement's noise.  The principle is to merge predictions from a  trusted model of
the dynamics of the  system with measurements, in order  to  1- efficiently filter the noise and 2-  get an accurate estimate of the (unknown) internal state of the system in real time.  The estimation relies on the computation at each time of a positive semi-definite matrix $P$,  that represents the covariance of the estimation error when the noises are assumed to be white and Gaussians. The differential equation verified by $P$ is a matrix Riccati equation. The Kalman filter is designed for linear dynamical systems, but when the system's trusted equations are non-linear,  the extended Kalman filter (EKF) is a natural extension of the Kalman filter based on a linearization of the equations around each estimate. In this case, $P$ verifies the same equation but with time-varying coefficients, and the computed positive semi-definite matrix is only an approximation of the true error covariance matrix. The stability analysis, both for the KF and EKF, generally relies on results on the Riccati flow on the set of positive semi-definite matrices, see \cite{wonham,song-grizzle-95,boutayeb}.

For non-linear systems, such as the one described by the aforementioned Riccati differential equation,  the contraction property may prove very useful. It essentially requires that the distance between two arbitrary  solutions does not expand. { When the  contraction is strict,   all solutions exponentially converge towards each other}.  Although a stronger mathematical property than global convergence {  to a particular solution}, as it implies  exponential convergence with a specified rate, it is sometimes an easier property to prove. In the control community it has been popularized by Lohmiler and Slotine \cite{slotine-auto}, { who} have presented it as  { a} useful tool to  analyze   some naturally stable systems, as in e.g. \cite{bonnabel-sepulchre-cdc11}, as well as for the design and stability analysis of non-linear  observers, see  \cite{slotine-auto,jouffroy-slotine,aghannan-rouchon-ieee03,bonnabel-geodesicIEEE}. { A potential interest of } the  contraction property is { that the use of a  non-Euclidean metric is often suggested by the dynamics and/or by the nonlinear nature of the state-space. }

A relevant example of non-Euclidean distance between covariance matrices is the natural distance of the cone {  of positive definite matrices}  \cite{faraut}. It enjoys several invariance properties, and appears in information geometry  as the Fisher distance between covariance matrices \cite{smith-2005}. P. Bougerol has proved in \cite{bougerol}  that the discrete-time Riccati equation is a strict contraction on the cone of positive definite matrices for the natural distance of the cone, see also \cite{Wojtkowski}. This generalizes to some extent the results of G. Birkhoff   on contraction properties  of positive maps on cones (see also \cite{bonnabel-sepulchre-cdc11} for observer design and contraction properties on cones). Those results have been partially extended in continuous time in \cite{jouffroy-slotine} using the same metric.  It seems to be the right non-Euclidean metric to analyse the stability and contraction properties of the Kalman filter.

When the dimension of the state space becomes very large, which is { e.g.} the case for discretized  systems described by partial differential equations, the Riccati equation of the KF and EKF become numerically intractable due to matrix storage and complexity of the updates. Essentially motivated by applications in weather forecasting and oceanography, where the covariance matrix can contain up to $10^7$ unknowns, the idea of projecting the state space vector onto a lower dimensional space has  appeared several decades ago e.g. \cite{Dee}. This operation reduces the storage and multiplication of matrices of size $n\times n$ to matrices of size $r\times n$ with $r\ll n$ going from quadratic complexity and memory requirement in the state space dimension $n$ to linear complexity in $n$. The low-dimensional subspace is supposed to capture the largest eigenvalues of the covariance matrix, where the error needs to be the most attenuated. { One such method is }  the so-called SEEK filter  \cite{pham,Rozier07}, which is a low-rank version of the Kalman filter (see also the recent paper \cite{hoteit}). { In contrast } to previous methods, the covariance matrix still has a fixed rank but its span is free to evolve at each step.  Unfortunately, the natural metric  is not defined on the boundary of the cone, that is, rank-deficient positive semi-definite matrices, and the hitherto analysis of the geometry of the Kalman Filter via its contraction properties for the natural  metric, does not readily apply to  low-rank Kalman filters.

In this paper, a novel  low-rank Kalman filter is proposed. The underlying Riccati flow evolves on the manifold of fixed rank symmetric positive semidefinite matrices. Exploiting the decomposition of any rank-$r$ positive semi-definite matrix $P=URU'$, where $U\in\RR^{n\times r}$ is a matrix whose columns form an orthonormal set, and $R$ is a $r-$dimensional positive  definite matrix,   a { Riemannian geometry}   recently introduced by the authors \cite{bonnabel-simax} has been shown to   { retain} many invariance properties of the natural metric of the cone.  As a consequence,  the low-rank Riccati flow admits contraction properties when the infinitesimal distances are measured with  this non-Euclidean  metric.  The contributions of this paper are threefold. First a continuous-time version of the contraction properties of the KF \cite{bougerol} for the natural metric of the cone is given. This contribution is rather tutorial, as the results of \cite{bougerol} in discrete time  imply  contraction in continuous time. This section can be skipped by a reader who would be familiar with contraction theory in continuous-time. The paper contains two novel contributions: a novel version of the low-rank Kalman filter defined in a proper geometric framework is proposed. And, the low-rank Kalman filter is proved to inherit some contractions properties of the Kalman filter for a recently introduced metric on the set of fixed-rank positive semi-definite matrices.   As a by-product stability properties are proved and  convergence speed rates are characterized around stable steady-states.

The paper is organized as follows. In Section 2, the contraction properties of the Kalman filter are recalled, and some results of \cite{bougerol} are shown to be provable directly in continuous time. In Section 3, a low-rank Kalman filter is introduced. In Section 4, we present the natural metric on the set of fixed-rank positive semi-definite matrices recently introduced by the authors in \cite{bonnabel-simax}. Finally, Section 5 addresses the contraction properties of the low-rank Kalman filter.

\subsection{Notation}
\begin{itemize}
    \item $\mathrm{P_+(n)}$ is the set of symmetric positive definite $n\times n$ matrices.
    \item $\mathrm{S^+(r,n)}$ is the set of symmetric positive semi-definite $n\times n$ matrices of rank $r\leq n$. We will only use this notation in the case $p<n$.
     \item $\mathrm{St(r,n)}=\mathrm{O(n)}/\mathrm{O(n-r)}$ is the  Stiefel manifold i.e. the set of $n\times r$
matrices with orthonormal columns: $U^TU=I_p$.
    \item span($A$) is the subspace of $\RR^n$ spanned by the columns of $A\in \RR^{n\times n}$.
    \item $T_X\mathcal M$ is the tangent space to the manifold $\mathcal M$ at $X$.
    \item $A'$ denotes the transpose of the matrix $A$.
\end{itemize}

\section{Contraction properties of the Kalman filter}
\label{contraction}

Consider the linear time-varying continuous-time system
\begin{equation}
\label{system}
\begin{array}{rcl}
\dotex x(t) & = & A(t) x(t) + G(t) w(t) \\
y(t) & = &  C(t) x(t) + H(t)\eta(t), \;
\end{array}
\end{equation} where $x(t) \in \RR^n$, $w(t) \in \RR^m$, $y(t),
~\eta(t) \in \RR^p$.  The random vectors $w(t)$ and $\eta(t)$ are
independent Gaussian white noise with zero mean and covariance
matrix equal to the identity. The matrices $A(t)$, $G(t)$, $C(t)$
and $H(t)$ have the appropriate dimensions and it is assumed that
$A(t)$ is invertible for all $t$.

The classical equations of the Kalman filter define a recursion for $$\hat x(t)=\mathbb E[x(t)|\{y(s)\}_{0\leq s\leq t}]$$the best estimate of the true state $x$, using the conditional error covariance matrix $P(t)$
\begin{align}
\dotex \hat x & =  (A-PC'(HH')^{-1}C) \hat x + PC'(HH')^{-1}y \label{KF}\\
\dotex P & =\Phi_t (P)= AP+PA'+GG'-PC'(HH')^{-1}CP\label{ricatti:eq}
\end{align}
The mapping $\Phi_t$ defines the continous matrix-valued Riccati differential Equation, and for each $P\in P_+(n)$,  $\Phi_t (P)$ is a tangent vector to $P_+(n)$ at $P$, i.e. $\Phi_t (P)\in T_PP_+(n)$.

\subsection{Natural metric of $P_+(n)$}
The geometry of the $n$-dimensional cone of symmetric positive
definite matrices  $P_+(n)$ has been well-studied in the literature.
The group $GL(n)$ acts transitively on this set via the following
action
\begin{equation}\label{gln:action} \gamma_A:P_+(n)\to
P_+(n),\quad P\mapsto APA' \end{equation} for any $A\in GL(n)$. If
$P$ is the covariance matrix of a gaussian vector of zero mean $x$,
then $\gamma_A(P)$ is the covariance matrix of the transformed
vector $Ax$. If $P,Q\in P_+(n)$ are two arbitrary points of the cone
there always exists $A\in GL(n)$ such that $Q=\gamma_A(P)$. This
property makes $P_+(n)$ a so-called homogeneous space under the
$GL(n)$ group action. The isotropy subgroup is the subgroup of $
GL(n)$ stabilizing the identity matrix i.e., $\{A\in GL(n),~
AIA'=I\}=O(n)$. As a general property of homogeneous spaces, the
following identification holds
$$
P_+(n)=GL(n)/O(n)
$$
The isotropy subgroup being compact, there exists a $GL(n)$-invariant Riemannian metric on $P_+(n)$ called the natural metric \cite{faraut}.  This metric is defined  as the usual scalar product  at the identity
$$
g_I(X_1,X_2) =\tr{X_1X_2}
$$where $X_1,X_2$ are tangent vectors at the identity, i.e. two symmetric matrices. If $Y_1,Y_2$ are two tangent vectors at $P\in P_+(n)$, the group action \eqref{gln:action} with $A=P^{-1/2}$ transports them to the tangent space at the identity, where $P^{-1/2}$ is defined as the symmetric square root of $P$.  The invariance of the metric then implies
\begin{align*}
g_{P}(Y_1,Y_2) &=g_I(P^{-1/2}Y_1P^{-1/2},P^{-1/2}Y_2P^{-1/2})=\tr{P^{-1}Y_1P^{-1}Y_2}\label{metric:def}
\end{align*}
The associated Riemannian distance is
$$
d_{P_+(n)}(P,Q)=\{\sum_1^n\log^2(\lambda_i)\}^{1/2}
$$
where $\lambda_1,\cdots,\lambda_n$ are the eigenvalues of $PQ^{-1}$. The main property of this metric is its invariance to conjugacy and inversion. For any $A\in GL(n)$ and $P,Q\in P_+(n)$
$$
d_{P_+(n)}(APA',AQA')=d_{P_+(n)}(P,Q)=d_{P_+(n)}(P^{-1},Q^{-1})
$$

Note that, this metric coincides with the Fisher Information Metric (FIM) for the following statistical inference problem: the available observations have a Gaussian
distribution with zero mean and a covariance matrix  parametrized by an unknown matrix in  $P\in P_+(n)$. The distance distorts the space to measure the amount of information between the distributions. For instance $d_{P_+(n)}(P,P+\delta P)\to \infty$ when $P\to 0$, for any $\delta P\succ 0$. This is understandable from an information point of view, as a Gaussian distribution with 0 mean and 0 covariance matrix, carries infinitely more information than a Gaussian with a strictly positive covariance matrix. Finally the FIM enjoys invariance properties to reparameterization i.e. $x\mapsto Ax$ in the case of 0 mean Gaussian distributions, and thus it is no surprise to see it coincide with the natural metric of the cone.

\subsection{Contraction property}

The notion of contraction \cite{slotine-auto} for a dynamical system described by the flow $\dotex x=f(x,t)$ can be interpreted as the (exponential) decrease of a the (geodesic) distance between two arbitrary points under  the flow.

\begin{defn}\label{contract:def}
Let $\dotex x=f(x,t)$ be a smooth dynamical system, defined on a
$C^1$ embedded manifold $M\in\RR$ equipped with a Riemanian metric
denoted by $g_x(v_1,v_2)$ on the tangent space at $x$. Let $X(x,t)$
denote the flow associated to $f$:
\begin{equation} X(x,0)=x,\quad
\dotex X(x,t)=f(X(x,t),t) \end{equation}
 Let $v$ be a tangent vector
at $x$. Let $\sigma_v(s),~0\leq s\leq 1$ be a geodesic emanating
from $x$, with direction $v$. Let $N$ be a subset of $M$. The system
is called a contraction for the metric $g$ on $N$ if for all $x\in
N$, and $v\in T_x M$,  we have \begin{equation}\label{contract:eq}
\dotex
\{g_x(\frac{dX(\sigma_v(s),t))}{ds},\frac{dX(\sigma_v(s),t))}{ds})\}\leq
-2\lambda
\{g_x(\frac{dX(\sigma_v(s),t))}{ds},\frac{dX(\sigma_v(s),t))}{ds})\}
\end{equation} where $\lambda\geq 0$. If $\lambda>0$ the contraction
is called strict.

\end{defn}

We have the following result \cite{slotine-auto,aghannan-rouchon-ieee03}:

\begin{thm}
Under the notation and assumptions of Definition \ref{contract:def}  let  $x_0,x_1\in N$. Assume there is a geodesic $\gamma(s),~0\leq s\leq 1$ linking $x_0$ and $x_1$, and for all $t\in[0,T]$ the transported geodesic $X(\gamma(s),t)\subset N$ for all $0\leq s\leq 1$, then if $d_g$ is the geodesic  distance associated to the metric we have
$$
d_g(X(x_0,t),X(x_1,t))\leq e^{-\lambda t}d_g(x_0,x_1),\quad\forall t\in[0,T]
$$which is an exponential decrease if $\lambda>0$.
\end{thm}
The proof relies on the fact that the  length of the curve $X(\gamma(s),t)$ transported by the flow, in the sense of the metric $g$, is the sum of all the length elements $\norm{{dX(\gamma(s),t)}/{ds}}_g$, and the infinitesimal property  \eqref{contract:eq}  can be extended to the length of the entire curve.

Strict contraction is a natural property to expect from a filter, as it implies the exponential forgetting of the initial condition. And it turns out that the discrete-time Riccati flow has been proved to be contracting for the natural metric of the cone $P_+(n)$ by Bougerol under standard assumptions, and it has also been proved to be a contracting for the same metric in continuous-time in \cite{jouffroy-slotine}, under more restrictive assumptions. Indeed if $\delta_P$ is a tangent vector at $P$ we have the following result
\begin{lem}
Assume that the process noise $GG'$ is stricly positive. Then the
Riccati equation \eqref{ricatti:eq} is a contraction in the sense of
Definition 1 for the natural metric of the cone. At each $P$ we have
$\lambda=\mu/p_{\text{max}}$ where $\mu$ is a lower bound on the
eigenvalues of $GG'$  and $p_{\text{max}}$ is an upper bound on the
eigenvalues of $P$.
\end{lem}
\paragraph{Proof}
We have
\begin{align*}
\dotex g_P(\delta P,\delta P)&=\dotex \tr{(P^{-1}\delta P)^2}\\&=2\tr{P^{-1}\delta P(\dotex P^{-1})\delta P}+2\tr{P^{-1}\delta P P^{-1}\dotex \delta P}
\end{align*}
where
$$
\dotex \delta P=A\delta P+\delta PA'-\delta PC'(HH')^{-1}CP-PC'(HH')^{-1}C\delta P
$$
and
$$\dotex P^{-1}=-P^{-1}\dot PP^{-1}=-P^{-1}
A-A'P^{-1}-P^{-1}GG'P^{-1}+C'(HH')^{-1}C$$
As a result
\begin{align*}
\dotex g_P(\delta P,\delta P)&=-2\tr{(P^{-1}\delta P  P^{-1}\delta P)GG'P^{-1}}\\&\qquad \quad-2\tr{(P^{-1}\delta P  P^{-1}\delta P)C'(HH')^{-1}CP}\\
&\leq -2\tr{(P^{-1}\delta P  P^{-1}\delta P)GG'P^{-1}}
\end{align*}
which proves that the Riccati flow is contracting in the sense of
\eqref{contract:eq} for the natural metric of the cone with
$\lambda=\mu/p_\mathrm{max}$ where $\mu$ is a lower bound on the
eigenvalues of $GG'$, that is assumed to be invertible, and
$p_\mathrm{max}$ is an upper bound on the eigenvalues of $P$.

The contraction property of the Riccati equation is in fact due to
the invariances enjoyed by the natural metric of the cone. This can
be understood the following way: the equation writes $\dotex P=
AP+PA'+GG'-PC'(HH')^{-1}CP$. The two first terms neither expand nor
contract as it is the differential form of the transformation
$\gamma_B(P)=(I+\tau A)P(I+\tau A)'$ which is an isometry for the
distance $d_{P_+(n)}$. The addition of a positive matrix $\dotex
P=Q$ where $Q$ is a positive matrix is neither expanding nor
contracting in the Euclidian space, but it contracts for the natural
metric as $g_P$ tends to dilate distances when $P$ becomes large.
Finally $\dotex P=-PC'(HH')^{-1}CP$ is a naturally contracting term
which is paramount in the theory of the Kalman filter and observers
(correction term). We have the following proposition, refining the
results of \cite{jouffroy-slotine} in the case of time-independent
coefficients:

\begin{prop}Suppose the matrices $A,C,G$ do not depend on the time $t$, $(A,C)$ is observable and  $G$ is full-rank. There exists a unique solution $Q\in P_+(n)$ of the algebraic Riccati equation defined by  $\Phi(Q)=0$ where $\Phi$ is given by \eqref{ricatti:eq}.  For any $R\in P_+(n)$, let $S_R$ be the ball center $Q$ and radius $d_{P_+(n)}(Q,R)$
$$
S_R=\{P\in P_+(n), ~d_{P_+(n)}(Q,P)\leq d_{P_+(n)}(Q,R)\}
$$
Let $M_R=\sup\{\norm{P}_2, P\in S_R\}<\infty$, and $\mu>0$ be the lowest eigenvalue of  $GG'$. Let $P_1(t),P_2(t)$ be two arbitrary solutions of the Riccati equation  $\dotex P=\Phi(P)$ initialized in $S_R$. Then  $P_1(t),P_2(t)\in S_R$ for any $t\geq 0$, and the following contraction result holds
$$
d_{P_+(n)}(P_1(t),P_2(t))=e^{-\frac{\mu }{M_R}t}d_{P_+(n)}(P_1(0),P_2(0))
$$
In particular $P(t)$ tends exponentially to the stationnary solution
$Q$ for the distance $d_{P_+(n)}$.
\end{prop}
\paragraph{Proof}The proof is a straightforward application of Theorem 1, where we used the fact that $d_{P_+(n)}$ is non-expanding and thus by Theorem 1  $P_1(t),P_2(t)$ remain in the compact set $S_R$.

\section{Low-rank Kalman filtering}

With ever growing size of measurable data and ever growing complexity of the models, the high computing cost and storage requirement of the Kalman filter (because of the need to implement the Ricatti equation) can become prohibitive. This has  especially been a problem in the field of  weather forecast and oceanography applications, where Kalman filters are a very natural and attractive tool, but where the models rely on discretized partial differential equations, and where the state vector as well as the measurement vector can become prohibitively large \cite{Rozier07}. Those problems resulted in modified versions of the extended Kalman filter where the state vector is projected on a low-dimensional space, see e.g.  \cite{Dee,pham}. Note that this filter could also meet applications in Extended Kalman Fitler based Simultaneous Localization and Mapping (e.g. \cite{dissanayake-2001}), where the Sparse Extended Information Filter, which is supposed to retain the largest terms in the covariance matrix, plays an important role.

In order to define  covariance matrices of  vectors that have been projected onto a subspace of dimension $r\ll n$, one  mainly thinks of two parameterization of positive semi-definite matrices of dimension $r$, i.e.  matrices of $S^+(r,n)$
$$
P=ZZ'=UR^2U'
$$where $P\in S^+(r,n)$, $Z\in \RR_*^{n\times r}$, $R\in \mathrm{P_+(r)}$, $U\in \mathrm{St({r,n})}$ (see figure). The use of the first decomposition goes back   to the early days of Kalman filtering, motivated
by numerical stability of the filter. Indeed, updating $Z$ and defining $P$ as $ZZ'$ enforces the constraint of a  positive symmetrical semi-definite covariance matrix $P$. The second decomposition $UR^2U'$ is  meaningful from a geometrical perspective and is suited to a simple statistical interpretation: $URU'$ is a flat ellipsoid, $U$ defines the subspace in which the ellipsoid lives, and $R$ defines the form of the ellipsoid. Such an ellipsoid is supposed to capture the principal directions of variability in the error covariance matrix and $R$ represents the error covariance matrix on the lower dimensional subspace. The set $S^+(r,n)$ has proved to be a Riemannian manifold, when equipped with appropriate metrics based on the first decomposition \cite{journee,meyer-11}, as well as the second  decomposition \cite{bonnabel-simax}. Another interesting metric  can also be found in \cite{vandereycken}.

\begin{figure*}
\includegraphics[width=1\textwidth]{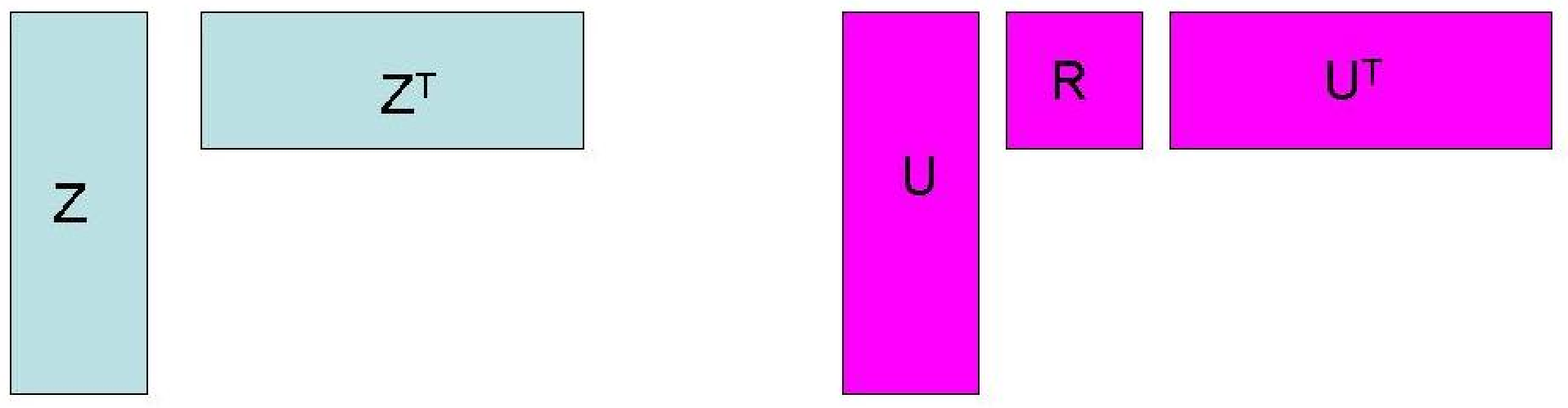}
\end{figure*}

In this paper, we seek to prove that the low-rank Riccati equation inherits some of the contraction properties of the full-rank Riccati equation. As the contraction property of the Riccati equation has been proved to stem from the invariance properties of the natural metric of the cone, we will focus on the second decomposition. Indeed,  metric \cite{bonnabel-simax}  enjoys several invariance properties and appears to be the most natural generalization of the natural metric of $P_+(n)$ to the set $S^+(r,n)$. { The matrix factorization}
\begin{align}\label{URU:eq}
P = U R U'
\end{align}
{ suggests to rewrite the Riccati equation as} \begin{equation}
\label{low-rankRicatti} \dotex U R U'=AU R U'+U R U'A'+GG'-U R
U'C'(HH')^{-1}CU R U'\end{equation} { However, this equation is not
rank-preserving because of the process noise term $GG'$.
Geometrically speaking, this means that $\dotex U R U'$ does not
belong to the tangent space  to the manifold $S^+(r,n)$ at $U R
U'$}. To circumvent this difficulty, \cite{pham} proposes to
{ disregard} the process noise, i.e.  $G=0$ so that the update is
rank-preserving. However, as seen in Section 2, process noise has
paramount stabilizing properties. This is why it is also proposed in
\cite{pham} to perturb a little the update with a forgetting factor,
in order to recover some stability properties. In this paper, we
propose a rank-preserving equation including process noise. To
preserve the rank,  the noise can be projected on the
low-dimensional space, i.e. $G\mapsto UG$. However, the invariance
properties discussed in the next section will in fact further
restrict the noise to be diagonal, i.e.  $G=\mu I$, leading to the
modified rank-preserving Riccati equation
\begin{align}\label{lr:ricatti:eq}
\dotex U R U^T=AU R U'+U R U'A'+\mu^2 UU'-U R U'C'(HH')^{-1}CU R U'
\end{align}
{ The contraction properties of  this low-rank Riccati equation are studied in the next sections.}

\section{An invariant metric on the manifold of fixed-rank positive semidefinite matrices}

From the previous sections, it seems natural to measure the contraction using the polar parametrization $P=URU'$.  However in this decomposition there is a catch, the underlying geometry is a quotient geometry because for any orthogonal matrix $O\in O(r)$ we have:

\begin{align}\label{trans:eq}
P = U R U'= (UO)(O'RO)(O'U')
\end{align}
The representation $P=URU'$ with $(U,R)\in \mathrm{St(r,n)}\times \mathrm{P_+(r)}$ is thus univoque up to the equivalence relation $(U,R)\equiv (UO,O'R' O)$ for any $O\in \mathrm{O(p)}$. The manifold $S^+(r,n)$ thus admits a quotient representation
$$S^+(r,n)\cong(\mathrm{St(r,n)}\times \mathrm{P_+(r)}/\mathrm{O(r)}$$
Let us introduce the metric proposed in \cite{bonnabel-simax}. Given a representative  $(U,R)$ of $P\in S^+(r,n)$, the tangent vectors of $T_PS^+(r,n)$ are given by the pair $(\Delta, D)$ where
\begin{equation}\label{tangent:space}
    \begin{aligned}\Delta&=U_\perp B,\quad\quad B\in \RR^{(n-r)\times r}\\D&=RD_0R\end{aligned}
\end{equation} such that $U_\perp\in St(r,n-r)$ , $U'U_\perp=0$, and $D_0\in \mathrm{Sym(p)}=T_I\mathrm{P_+(r)}$. The chosen metric of $S^+(r,n)$ is merely  the sum of the infinitesimal distance between subspaces and  between low-rank positive definite matrices of $\mathrm{P_+(r)}$:
\begin{align}\label{metric:def} g_{(U,R^2)}((\Delta_1,D_1),(\Delta_2,D_2))=\tr{\Delta_1^T\Delta_2}+\tr{R^{-1}D_1R^{-2}D_2R^{-1}},
\end{align}generalizing the natural metric of the cone in a natural way. According to   \cite{bonnabel-simax},
the space $S^+(r,n)\cong(\mathrm{St(r,n)}\times \mathrm{P_+(r)}/\mathrm{O(r)}$ endowed with the metric \eqref{metric:def} is a Riemannian manifold with horizontal space $$\mathcal
H_{(U,R^2)}=\{(\Delta,D): \Delta=U_\perp B, ~B\in \RR^{(n-r)\times r},
~D=RD_0R,~D_0\in \mathrm{Sym(p)}\}$$
Furthermore, the metric proved to inherit some invariance properties of the natural metric, namely invariance to orthogonal transformations, dilations, and pseudo-inversion.

\section{Contraction properties of the low-rank Kalman filter}

The Riccati low-rank flow \eqref{lr:ricatti:eq} defines a tangent vector to the manifold. It corresponds to the following flow defined in the horizontal space
\begin{align}
\dotex U& =(I-UU')AU \label{powerflow:eq}\\
 \dotex R&=A_UR+RA_U' +\mu^2 I-R C_{U}'(HH')^{-1}C_U R \label{ricatti2:eq}
\end{align}
where $A_U=U'AU$, $C_U=CU$ are the matrices projected in span $U$.

The system has a triangular structure, that is,
  the first equation is independent from the second one. The dynamics of $U$ are known as the (non symmetric) Oja flow.  In turn, $R$ verifies a time-varying Riccati equation on the lower dimensional cone $P_+(r)$.

Topological obstructions seem to prevent us from proving the low-rank Kalman filter is a global contraction on $\mathrm{S^+(r,n)}$. Indeed, even in the simplest case $r=1$ and $n=2$, the flow  \eqref{powerflow:eq} defines a vector field that can not be a contraction on  the whole circle  St$(1,2)$. This would imply there exists a globally converging vector field on the circle which is impossible. Thus, global contraction seems to be out of reach, and we are going to prove instead local stability and contraction around any dominant invariant  subspace.

\subsection{Local contraction of subspace flow and low-rank Riccati flow}

The flow \eqref{powerflow:eq} has been well studied under the assumption $A=A'$. In this case, the flow is the gradient flow of the generalized Rayleigh quotient cost function  $r_A(U)= \tr{U'A U}$ for the natural metric on the Stiefeld manifold $St(r,n)$. Because the cost is invariant to rotations, the gradient flow can also be interpreted as a gradient flow on the Grassmann manifold $\mathrm{Gr(r,n)} \approx \mathrm{St(r,n)}/\mathrm{O(r)}$, i.e. the manifold of $r$-dimensional subspaces in a $n$-dimensional ambiant space. This subspace flow has been studied by Oja \cite{oja-92} and by Wei-Yong, Helmke and Moore \cite{helmke}. If the eigenvalues of $A$ are distinct, the flow converges for almost every initial condition to the dominant subspace of $A$, i.e. the eigenspace associated with the $r$ largest eigenvalues of $A$.

For any matrix $A$, {\it local} exponential convergence of the dominant eigenspace is ensured provided that there is a gap between the $r$ first eigenvalues of the symmetric part of $A$, and the others, as shown by the following lemma.

\begin{lem}\label{grass:contr}
The matrix $(A+A')/2$ is real symmetric and can be diagonalized. Let $\lambda_U$ be the gap between the $r$-th and the $r+1$-th eigenvalues.  Let $U_r\in \mathrm{St(r,n)}$ represent an orthonormal basis of the the subspace spanning the $r$ dominant eigenvalues of the symmetric part of $A$. If $\lambda_U>0$, $U_r$ is an exponentially stable equilibrium of the subspace flow \eqref{powerflow:eq} in the sense of the Grassman  metric.
\end{lem}

\paragraph{Proof}
 Let $U_\perp$ span the orthocomplement of $U_r$ and $\delta U=U_\perp K$
 with $K\in \RR^{(n-r)\times r}$.  $U_r$ is an
 equilibrium, and the first-order variation of \eqref{powerflow:eq} reads$$
\dotex \delta U=(I-U_rU_r')A\delta U-\delta UU_r'AU_r-U_r\delta U'AU_r
$$
which implies
\begin{align*}
\dotex\tr{\delta U'\delta U}&=2\dotex \tr{\delta U'\dotex \delta U}\\
&=2\tr{\delta U' A\delta U- \delta U'\delta UU_r'AU_r}\\
&=\tr{\delta U' (A+A')\delta U- \delta U'\delta UU_r'(A+A')U_r}\\
&=\tr{K'U_\perp (A+A'U_\perp K- K'KU_r'(A+A'U_r}
\\&\leq -2\lambda_U\tr{KK'}\leq-2\lambda_U\tr{\delta U'\delta U}
\end{align*}

We thus see that every solution that converges to the dominant eigenspace of the symmetric part of $A$, is eventually contracting for the natural metric of the Grassmann (subspace) manifold. We moreover naturally have the following result.

\begin{lem}\label{co:prop}The  dynamics \eqref{powerflow:eq}-\eqref{ricatti2:eq} is a strict contraction for the metric \eqref{metric:def} on any compact subset of $S^+(r,n)$ on which the span  is held fixed. 
\end{lem}

\subsection{Eventual contraction in the case of symmetric $A$}\label{eventual:subsec}

As mentionned above, system \eqref{powerflow:eq}-\eqref{ricatti2:eq} has a triangular structure. The dynamics of $U$ is independent of the dynamics of $R$, and  under the assumption that $A=A'$ and the eigenvalues of $A$ are distinct, the flow converges for almost every initial condition to the dominant subspace of $A$.  Once $U(t)$ has converged to $U_\infty=U_r$, the matrices $A_U$ and $C_U$ are projections of $A$ and $C$ on a fixed subspace $U_r$. As a result, if $A$ and $C$ are time-independent,
 \eqref{ricatti2:eq} asymptotically becomes a Riccati equation with time-independent coefficients, implying convergence of $R(t)$ from Proposition 1, as soon as $(A_{U_\infty},C_{U_\infty})$ is observable. Those two latter facts allow to conclude the cascaded system \eqref{powerflow:eq}-\eqref{ricatti2:eq}  converges to fixed values $(U_\infty,R_\infty)$, as   explosions in finite time can not occur (it is easily checked that the trace of $R$ can grow at most exponentially).

As it has been shown that the subspace flow is a contraction around the final value $U_\infty$, and when $U$ is  held equal to $U_\infty$ the low-rank  Riccati equation is a contraction, we can thus conclude the low-rank Riccati flow ``eventually" possesses some contraction properties under the above assumptions on $A$.

\begin{rem}The complete filter  also involves the equation for the estimated mean  $\hat x$ of $x$, that writes $\dotex \hat x  =  (A-KC) \hat x + Ky$ with $K=URU'C'(HH')^{-1}$. Its convergence properties are a consequence of the contraction properties of the low-rank Riccati equation. Indeed, as the gain $K$ asymptotically converges, it can be proved the mean error $\hat x-x$ projected on the final subspace spanned by $U_\infty$ converges to zero under the additional assumption that the output map satisfies $CU_\infty U_\infty'=C$.
\end{rem}

\subsection{A counter-example for skew-symmetric $A$}
Convergence properties are necessarily weaker in the case of a nonsymmetric matrix $A$, as illustrated by the following example.  Consider \eqref{powerflow:eq} on St(1,3) with $A$ skew-symmetric. Then $U=x\in\RR^3$ is a vector on the unit sphere and it satisfies
$$
\dotex x=\omega\wedge x
$$for some $\omega\in\RR^3$. It is clear that this flow does not possess
any contraction property as an initial angle between two vectors remains unchanged by the
flow, and thus the distance in the Stiefel manifold is constant over time.

For $A\in\RR^{n\times n}$ such that $A-A'\neq 0$, the skew symmetric part of A induces a rotation
of the subspace spanned by U, and one should not  expect neither
convergence, nor contraction. It seems reasonable to hope for
ultimate contraction, even with general A, as soon as the associated full-rank equation converges to a stationnary state.
Proving this statement is beyond the scope of the present paper.

\subsection{Proposed implementation in discrete-time}

In order to meet the implementation constraints, the proposed filter must be written in discrete time. The transposition is not straightforward, as the rank constraint must be preserved numerically at each step, i.e. the updated covariance matrix must remain on the manifold. It means the updated $U$ must have orthonormal columns, and the updated $R$ must be positive definite. The so-called retractions (e.g., \cite{absil-book}) conveniently map the tangent space at any point, to the manifold. Using the retractions adopted in the work  \cite{absil-book,meyer-11}, and the usual discrete Riccati equation we propose the following implementation in discrete time:
\begin{align*}
 U^+& =\text{qf}(\tilde AU-dt (UU'AU)) \\
 R^+&=\tilde A_UR\tilde A_U'-\tilde A_URC_U'(C_URC_U'+dtHH')^{-1}C_UR\tilde A_U'+dt\mu^2I
 \end{align*}
where $\tilde A=A+dtI$  and qf() extracts the orthogonal factor in the QR decomposition of its argument. Note this latter operation boils down to orthonormalizing the factor $U$ at each step.

\section{Conclusion}

In this paper we have analyzed contraction properties of the low-rank Kalman filter via a recently introduced metric that extends the natural metric (or Fisher information metric) to the facets of the cone of positive definite matrices. In the process, the low-rank Kalman filter equations were a little modified and put on firm geometrical ground.  Here are some  concluding remarks and challenges that  could deserve attention:
\begin{itemize}
\item In discrete-time \cite{bougerol}, the contraction rate of the Riccati flow for the natural metric of $\mathrm{P_+(n)}$ only depends on the coefficients $A,B,C,D$. However, in continuous-time it depends on the point $P\in \mathrm{P_+(n)}$ at which it is evaluated  (see Lemma 1). This difference is a little surprising and would deserve to be better understood. Note also, that if one can prove simply there is an arbitrary large invariant set for the flow under detectability and observability conditions, Lemma 1 implies the contraction rate is uniform and hence convergence by the fixed point theorem.
\item Thus, one could hope to prove stronger results on the contraction properties of the  low-rank Riccati flow in discrete-time. However, computing explicitely the distance between two low-rank matrices with the metric proposed remains an open problem. \cite{bonnabel-simax} proposes special curves in $\mathrm{S^+(r,n)}$ that approximate the geodesics. Proving stronger results on the contraction  of thoses special curves under the Riccati flow seems to be straightforward enough.
\item In this paper we proved  contraction properties of the flow on a  subset of the manifold. Although global contraction seems out of reach because of topological obstructions, some ``eventual"  contraction properties were underlined in the special case where $A$ is symmetric. Without this assumption, Oja flow, and thus the low-rank Riccati flow, do not necessarily converge.
\item The filter proposed in this paper requires a reduced number of numerical operations and
storage capacity. It would be of interest to test its efficiency in a particular large-scale
application and to evaluate its relative merits with respect to alternative approaches
that focus on sparsity of the covariance matrix (e.g. \cite{Thrun02simultaneousmapping}).
\end{itemize}

\bibliographystyle{plain}

\begin{thebibliography}{10}

\bibitem{absil-book}
P.A. Absil, R.~Mahony, and R.~Sepulchre.
\newblock {\em {Optimization Algorithms on Matrix Manifolds}}.
\newblock {Princeton University Press}, 2007.

\bibitem{aghannan-rouchon-ieee03}
N.~Aghannan and P.~Rouchon.
\newblock An intrinsic observer for a class of lagrangian systems.
\newblock {\em IEEE Trans. on Automatic Control}, 48(6):936--945, 2003.

\bibitem{vandereycken}
Stefan~Vandewalle Bart~Vandereycken, P.-A.~Absil.
\newblock A riemannian geometry with complete geodesics for the set of positive
  semidefinite matrices of fixed rank.
\newblock {\em Submitted}.

\bibitem{bonnabel-geodesicIEEE}
S.~Bonnabel.
\newblock A simple intrinsic reduced-observer for geodesic flow.
\newblock {\em IEEE Trans. on Automatic Control}, 55(9):2186 -- 2191, 2010.

\bibitem{bonnabel-sepulchre-cdc11}
S.~Bonnabel, A.~Astolfi, and R.~Sepulchre.
\newblock Contraction and observer design on cones.
\newblock In {\em IEEE Conference on Decision and Control 2011 (CDC11)}, 2011.

\bibitem{bonnabel-simax}
S.~Bonnabel and R.~Sepulchre.
\newblock Riemannian metric and geometric mean for positive semidefinite
  matrices of fixed rank.
\newblock {\em SIAM J. matrix anal. appl.}, 31, 2009.

\bibitem{bougerol}
P.~Bougerol.
\newblock Kalman filtering with random coefficients and contractions.
\newblock {\em Siam J. Control and Optimization}, 31:942--959, 1993.

\bibitem{boutayeb}
M.~Boutayeb, H.~Rafaralahy, and M.~Darouach.
\newblock Convergence analysis of the extended kalman filter used as an
  observer for nonlinear deterministic discrete-time systems.
\newblock {\em IEEE transactions on automatic control}, 42, 1997.

\bibitem{Dee}
D.P. Dee.
\newblock Simplification of the kalman filter for meteorological data
  assimilation.
\newblock {\em Quart. J. Roy. Meteor. Soc.}, 10:365--384, 1990.

\bibitem{dissanayake-2001}
G.~Dissanayake, P.~Newman, H.F. Durrant-Whyte, S.~Clark, and M.~Csobra.
\newblock A solution to the simultaneous localisation and mapping (slam)
  problem.
\newblock {\em IEEE Trans. Robot. Automat.}, 17:229--241, 2001.

\bibitem{hoteit}
M.~El~Gharamti, I.~Hoteit, and S.~Sun.
\newblock Low-rank kalman filtering for efficient state estimation of
  subsurface advective contaminant transport models.
\newblock {\em Journal of Environmental Engineering}, 2011.

\bibitem{faraut}
J.~Faraut and A.~Koranyi.
\newblock {\em Analysis on Symmetric Cones}.
\newblock Oxford Univ. Press, London, U.K., 1994.

\bibitem{jouffroy-slotine}
J.~Jouffroy and J.J.E. Slotine.
\newblock Methodological remarks on contraction theory.
\newblock In {\em 43rd IEEE Conference on Decision and Control}, 2004.

\bibitem{journee}
M.~Journee, P.-A.~Absil F.~Bach, and R.~Sepulchre.
\newblock Low-rank optimization on the cone of positive semidefinite matrices.
\newblock {\em SIAM Journal on Optimization}, 20(5):2327--2351, 2010.

\bibitem{Kalman-1961}
R.~Kalman and R.~Bucy.
\newblock New results in linear filtering and prediction theory.
\newblock {\em Basic Eng., Trans. ASME, Ser. D,}, 83(3):95--108, 1961.

\bibitem{slotine-auto}
W.~Lohmiler and J.J.E. Slotine.
\newblock On metric analysis and observers for nonlinear systems.
\newblock {\em Automatica}, 34(6):683--696, 1998.

\bibitem{meyer-11}
G.~Meyer, S.~Bonnabel, and R.~Sepulchre.
\newblock Regression on fixed-rank positive semidefinite matrices: a riemannian
  approach.
\newblock {\em Journal of Machine Learning Reasearch (JMLR)}, 12:593--625,
  2011.

\bibitem{oja-92}
E.~Oja.
\newblock Principal components, minor components, and linear neural networks.
\newblock {\em Neural Networks}, 5:927 -- 935, 1992.

\bibitem{pham}
D.T. Pham, J.~Verron, and M.C. Roubaud.
\newblock A singular evolutive extended kalman filter for data assimilation in
  oceanography.
\newblock {\em Journal of Marine Systems}, 16:323--340, 1998.

\bibitem{Rozier07}
D.~Rozier, F.~Birol, E.~Cosme, P.~Brasseur, J.-M. Brankart, and J.~Verron.
\newblock A reduced-order kalman filter for data assimilation in physical
  oceanography.
\newblock {\em SIAM Rev.}, 49:449--465, 2007.

\bibitem{smith-2005}
S.T. Smith.
\newblock Covariance, subspace, and intrinsic cramer-rao bounds.
\newblock {\em IEEE-Transactions on Signal Processing}, 53(5):1610--1629, 2005.

\bibitem{song-grizzle-95}
Y.K. Song and J.W. Grizzle.
\newblock The extended kalman filter as a local asymptotic observer.
\newblock {\em Estimation and Control}, 5:59--78, 1995.

\bibitem{Thrun02simultaneousmapping}
Sebastian Thrun.
\newblock Simultaneous mapping and localization with sparse extended
  information filters: Theory and initial results.
\newblock 2002.

\bibitem{helmke}
Y.~Wei-Yong, U.~Helmke, and J.B. Moore.
\newblock Global analysis of oja's flow for neural networks.
\newblock {\em IEEE Trans. on Neural Networks}, 5:674--683, 1994.

\bibitem{Wojtkowski}
M.P. Wojtkowski.
\newblock Geometry of kalman filters.
\newblock {\em Journal of geometry and symmetry in physics}, 2007.

\bibitem{wonham}
W.M. Wonham.
\newblock On a matrix ricatti equation of stochastic control.
\newblock {\em SIAM J. Control}, 6, 1968.

\end{thebibliography}

\end{document}